\documentclass{article}
\usepackage{amssymb}
\usepackage{lineno}
\usepackage{microtype}
\usepackage{graphicx}
\usepackage{subfigure}
\usepackage{booktabs} 
\usepackage{amsmath, amsthm, amssymb, multirow,fullpage}
\usepackage{bm}

\usepackage[normalem]{ulem}
\RequirePackage[authoryear,round]{natbib}
\usepackage[dvipsnames,table,xcdraw]{xcolor}
\usepackage[colorlinks=true, allcolors=blue,backref=page]{hyperref}
\usepackage{setspace}
\usepackage{epsfig,psfrag,url,verbatim,mathrsfs}
\usepackage{graphics}
\usepackage{graphicx}
\usepackage{extarrows}
\usepackage{ctable}
\usepackage{multirow}
\usepackage{tabu}
\usepackage{amsfonts}
\usepackage{algorithm}
\usepackage{algpseudocode}
\usepackage{pifont}

\usepackage{wrapfig}
\usepackage[colorlinks=true, allcolors=blue,backref=page]{hyperref}
\renewcommand*{\backref}[1]{}
\setcitestyle{citesep={,}}
\usepackage{ctable}
\usepackage[parfill]{parskip}
\usepackage[nameinlink]{cleveref}

\numberwithin{equation}{section}

\newtheorem{definition}{Definition}[section]
\newtheorem{example}{Example}[section]
\newtheorem{proposition}{Proposition}[section]

\newtheorem{theorem}[proposition]{Theorem}
\newtheorem{lemma}[proposition]{Lemma}

\newtheorem{remark}[proposition]{Remark}

\allowdisplaybreaks

\crefname{assumption}{Assumption}{Assumptions}
\crefname{lemma}{Lemma}{Lemmata}
\crefname{definition}{Definition}{Definitions}
\crefname{example}{Example}{Examples}
\crefname{proposition}{Proposition}{Propositions}
\crefname{theorem}{Theorem}{Theorems}
\crefname{corollary}{Corollary}{Corollaries}
\crefname{claim}{Claim}{Claims}
\crefname{algorithm}{Algorithm}{Algorithms}
\crefname{Def}{Definition}{Definition}
\crefname{remark}{Remark}{Remarks}
\crefname{figure}{Figure}{Figures}
\crefname{section}{Section}{Sections}
\crefname{equation}{}{}
\crefname{table}{Table}{Tables}

\usepackage[noblocks,auth-sc]{authblk}
\usepackage[normalem]{ulem}

\usepackage{accents}
\renewcommand{\bar}{\overline}
\renewcommand{\tilde}{\widetilde}
\renewcommand{\hat}{\widehat}

\renewcommand{\gcd}{\operatorname{GCD}}
\newcommand{\supp}{\operatorname{supp}}

\newcommand{\lcm}{\operatorname{LCM}}
\newcommand{\dist}{\operatorname{dist}}
\newcommand{\balldist}{\mathbb{B}_{\dist}}

\newcommand{\R}{\mathbb{R}}
\newcommand{\Z}{\mathbb{Z}}
\newcommand{\calD}{\mathcal{D}}

\newcommand{\calU}{\mathcal{U}}
\newcommand{\calX}{\mathcal{X}}
\newcommand{\basea}{\calX_{\alpha}}
\newcommand{\aggFun}{\bm{\vartheta}}
\newcommand{\ceil}[1]{\left\lceil #1 \right \rceil }

\newcommand{\bmu}{\bm{u}}
\newcommand{\bmx}{\bm{x}}
\newcommand{\bmy}{\bm{y}}

\newcommand{\utility}{u}
\newcommand{\utilvec}{\bmu}
\newcommand{\unfair}{\phi}

\title{A Framework for Fair Decision-making Over Time with Time-invariant Utilities }
\author[1]{Andrea Lodi}
\author[2]{Sriram Sankaranarayanan}
\author[3]{Guanyi Wang}
\affil[1]{Cornell Tech and Technion -- IIT}
\affil[2]{Indian Institute of Management, Ahmedabad}
\affil[3]{National University of Singapore}
\date{}

\begin{document}

\maketitle
\begin{abstract}
Fairness is a major concern in contemporary decision problems. 
In these situations, the objective is to maximize fairness while preserving the efficacy of the underlying decision-making problem. 
This paper examines repeated decisions on problems involving multiple stakeholders and a central decision maker. 
Repetition of the decision-making provides additional opportunities to promote fairness while increasing the complexity from symmetry to finding solutions. This paper presents a general mathematical programming framework for the proposed fairness-over-time (FOT) decision-making problem. 
The framework includes a natural abstraction of how a stakeholder's acquired utilities can be aggregated over time. 
In contrast with a natural, descriptive formulation, we demonstrate that if the aggregation function possesses certain basic properties, a strong reformulation can be written to remove symmetry from the problem, making it amenable to branch-and-cut solvers. 
Finally, we propose a particular relaxation of this reformulation that can assist in the construction of high-quality approximate solutions to the original problem and can be solved using simultaneous row and column generation techniques.
\end{abstract}

\section{Introduction}
Choosing a good alternative from a given (discrete or continuous) set of alternatives has been the central concept underlying numerous decision-making problems.
Typically, the best alternative is the one that maximizes or minimizes an objective that reflects the decision-maker's interests.
In many instances, however, there may be a single \emph{decision-maker} but multiple \emph{stakeholders} of the decision.
In such situations, it may be crucial for the decision maker to select not only an option of optimal objective for the whole society but also one that is {\em fair} and equitable to each party involved.
This spawned a new field of study on fairness.
\paragraph{Fairness in general decision-making problems}
The majority of papers on fairness focus primarily on \emph{resource-allocation problems} \citep{bai2022fair,flanigan2021fair,brandt2016handbook,aziz2020developments,bampis2018fair}.
These are problems involving the allocation of resources to a set of stakeholders, who could have heterogeneous resource valuations and non-additive resource bundle valuations.
In such situations, protocols are developed to enhance various notions of fairness in the proposed allocation.

However, not all decision-making problems involving fairness considerations are resource allocation problems. 
In typical resource allocation problems, given a set of items, one has to identify the destination where each item goes. A key point here is that each stakeholder's utility depends only upon their own basket of goods. On the other hand, consider the problem of positioning an ambulance vehicle in a town. The same vehicle can provide a certain utility to a stakeholder located 100 meters from the vehicle, some other utility to another stakeholder located 150 meters from the vehicle, and so on.
Moreover, unlike resource allocation problems, there is no single stakeholder to whom the ambulance vehicle is assigned. However, the position of the vehicle gives different utilities to different stakeholders.

Consider the following circumstance as another example. A pick-up vehicle must depart from a hub, pick up $n$ stakeholders from their respective locations, and then return with everyone to the hub, i.e., performing a trip of the well-known traveling salesman problem (TSP). 
The disutility of each stakeholder could depend on how long they do spend in the pick-up vehicle before being dropped off at their hub. Again, no resource is allocated to any stakeholder here. 
However, the path taken to pick up different stakeholders and return to the hub determines the disutility obtained by each of them.

Although these novel problems do not fall into the frameworks examined in the literature on fair resource allocation, it is straightforward to model them using a mathematical programming framework. 
Therefore, a fair decision-making approach, based on mathematical programming is warranted for problems requiring greater modeling complexity.

\paragraph{Fairness over time}
While one-shot decision-making problems are intriguing, usually, little can be 
done to remedy the unfairness. 
For example, in the earlier-mentioned TSP example, if the trip must be made only once, there will always be some stakeholder picked up first and must thus visit all the others before reaching the hub. 
As a result, there is little that can be done to improve fairness. 
However, if such trips are repeated, we notice two things. 
First, using the same cost-minimizing trip in each decision epoch (or period) is even more unfair to the stakeholder picked up first. 
Second, in another decision epoch, we can choose a slightly less cost-effective trip in which a different stakeholder is picked up first. 
Thus, over multiple rounds of decision making, stakeholders could take turns being selected first, thereby reducing unfairness. 
Moreover, it is easy to include any cost restrictions as a constraint in such models.

Besides, \citet{lodi2021fairness} provide explicit mathematical examples where \emph{perfect} fairness (based on a selected definition of fairness) is not possible when there is only one round of decision making. However, given multiple rounds of decision making, perfect fairness is possible. 

\paragraph{Contributions} 
In this paper, we present a mathematical programming framework applicable to a wide class of {\em repeated} decision-making problems where fairness is of importance. 
This generalizes the literature on resource allocation problems by permitting general decision-making issues. 
This also expands the work in \citet{lodi2021fairness} by allowing for more general \emph{utility aggregation} by each stakeholder. 
In particular, \citet{lodi2021fairness} assume that if a stakeholder gets different utilities over $T$ decision epochs, then the stakeholder always considers the {\em average} of those utilities as the representative utility, or what we call the {\em aggregated} utility. 
Motivating that aggregating by averaging might not be the most appropriate aggregation in some cases (\cref{ex:whyGenAgg}), we generalize and provide an abstraction for the class of aggregation functions that can  be useful.
More importantly, we provide efficient integer programming formulations for solving such problems. 
Given that such an integer program will necessarily be symmetric across each decision period, we provide an alternative formulation that reduces symmetry. 
This guarantees computational speed-ups for solving the problem.

Next, we show that one can solve a useful relaxation of the problem and get deep insights into the original problem. 
While it is known that solving a relaxed version of a problem provides useful bounds, we show that even stronger insights are possible for the case in hand.

{In particular, we prove that solving a relaxed version of the problem answers the following questions. 
\begin{enumerate}
\item What is the maximum fairness (strictly, minimum unfairness, as we will work using \emph{unfairness} -- or disutility -- functions) one can obtain, given unlimited number of decision periods?
\item \label{en:2} What is the maximum number of decision periods required to obtain an approximately fair solution? More precisely, let $L$ be the largest Lipschitz constant associated with any of the fairness functions involved in the decision variables for the underlying decision-making problem. Let $\hat{\unfair}$ be the largest fairness (smallest unfairness) possible allowing unlimited rounds of decision making. We show that, in $T = \ceil{\frac{L}{\epsilon}}$ rounds, it is possible to achieve unfairness less than $\hat{\unfair}+\epsilon$, for any $\epsilon > 0$.
\end{enumerate}
In particular, to answer question \ref{en:2}, we provide a constructive answer, i.e., a sequence of decisions for $T$ decision periods, so that the unfairness at the end of $T$ epochs is at most $\hat{\unfair}+\epsilon$. 
}
{

{
\paragraph{Literature Review} 
Fair resource allocation  is the problem of dividing a resource or a collection of resources among agents who derive value out of those resources.
They are broadly classified as {\em divisible} and {\em indivisible} resources. Divisible resource allocation corresponds to the so-called {\em cake-cutting protocols}. 
In this setting, a heterogeneous cake that can be cut anywhere has to be {\em fairly} divided among stakeholders with dissimilar valuations for different parts of the cake. 
In contrast, indivisible resource allocation corresponds to allocating a collection of finitely many items to different stakeholders. 
This problem is studied in various settings. For example, one can restrict oneself to the cases where each stakeholder associates a value for each item. Then, the total value of a collection of items is the sum of the values of the individual items. Generally, one can define a valuation function from each subset of the total pool of objects to reals. We refer the reader to \citet{brandt2016handbook} for a comprehensive review of these topics.

Results proving the existence of {\em fair} allocations and our ability to find them fast are shown for various definitions of fairness and various assumptions on valuations. These definitions of fairness range from envy freeness, envy freeability, envy freeness up to one good, and proportionality among others. 
Negative results indicating NP-hardness of deciding the existence of fair allocations under certain combinations of fairness definitions and valuation assumptions are not uncommon either. 

One of the primary concerns in this area is what is the loss of social welfare (i.e., efficiency) as one requires fairness. \citet{bertsimas2011price} answer this question, providing tight bounds for various notions of fairness when social welfare is defined as the sum of the utilities obtained by individual stakeholders. 
However, they restrict themselves to a single decision period. 

\citet{lodi2021fairness} consider decision making over multiple periods for a more general problem.
In particular, they consider the problem of allocating ambulance vehicles to different stations or bases. 
Stakeholders located in different parts of the cities achieve different utilities due to the presence of an ambulance at a base, depending upon their distance to the base. 
For this problem, they provide an integer programming algorithm to identify a sequence of allocations such that fairness is enhanced over multiple rounds of decision making. As stated earlier, the paper assumes that the aggregation of utilities over multiple decision periods is obtained by averaging them.

\citet{salem2022enabling} extends the work in \citet{lodi2021fairness} by considering an online version of the problem.
Like before, fairness is enhanced over multiple decision periods. 
However, the exact set of utility functions in any period is not known \emph{a priori}. 
While restricting themselves to only convex feasible sets and convex objectives, they provide algorithms to identify fair solutions, which they prove to be not much worse than a static decision taken with the benefit of hindsight.

}

Moreover, \citet{delage2019dice,delage2022value,castiglioni2022unifying} independently explore the insights/connections between randomization and convex relaxation.
    \citet{delage2019dice} investigate circumstances under which random decisions strictly outperform deterministic decisions. They show that several convex risk measures, such as mean (semi-)deviation and mean (semi-)moment measures, can benefit from randomization if the set of feasible decisions is nonconvex. Moreover, they show that under a mild continuity property, for any ambiguity-averse risk measure, one can find a decision problem with a nonconvex feasible set in which a randomized strategy strictly dominates all pure strategies. 
    \citet{delage2022value}  study the value of randomization in mixed-integer \emph{distributionally robust} optimization problems. 
    They provide improvement bounds based on continuous relaxation and identify the conditions under which these bounds are tight. On the algorithmic side, they propose a finitely convergent, two-layer, column-generation algorithm that iterates between identifying feasible solutions and finding extreme realizations of the uncertain parameter. 
    \citet{castiglioni2022unifying} utilize the so-called strategy mixtures to study online learning problems by introducing the set of probability measures on the original decision set. In that paper, the connection between randomization and convexification has also been captured in the analysis of strong duality, although it has not been explicitly presented. 

\paragraph{Paper Organization} We organize the paper as follows. \cref{sec:Motivation} motivates the problem of interest, and provides abstract definitions of unfairness functions and aggregation functions. With these definitions, we define the abstract fairness-over-time problem.
Next, in \cref{sec:PER}, we define the probabilistic-equivalent reformulation, which reformulates the problem of interest in a form that counters the symmetry in the original formulation.
\cref{sec:ProbRel}, then defines a certain relaxation of the probabilistic-equivalent reformulation and shows the useful information obtained on solving the relaxed version of the problems. 
\cref{sec:PER,sec:ProbRel} together contains most of the important results of the paper.
Finally, \cref{sec:Disc} concludes the paper by discussing the results, limitations, and future directions.

\section{Motivation \& Problem Setting} \label{sec:Motivation}

{
Consider the classical version of a decision-making optimization problem. 
We refer to the decision-maker's preferences, excluding the fairness requirements as {\em efficiency}. 
Thus, the objective that is maximized in the traditional sense is referred to in this manuscript as the {\em efficiency metric.}
} 
Such a problem can be abstractly written in the form
\begin{align}
    \max_{\bmx} ~ c(\bmx) ~~\text{s.t.}~~ \bmx \in \calX, \label{eq:dm-opt}
\end{align}
where $\calX$ denotes the feasible set of the decisions $\bmx$, and $c: \calX \to \mathbb{R}$ denotes the efficiency metric. 

{
We assume that there are $n\in \mathbb N$ stakeholders who derive heterogeneous utility out of each decision $\bmx\in\calX$. 
We formalize this notion by defining the utility function below. 
\begin{definition}[Utility Function \citep{lodi2021fairness}]
For each $i \in [n]$, let $\utility_i :\calX\to \calU_i \subseteq \mathbb{R}$ be the utility function associated with any decision for stakeholder (client) $i$. 
We use $\utilvec(\bmx) = (\utility_1(\bmx), \ldots, \utility_n(\bmx))^{\top}$ to denote the vector of utilities associated with the decision $\bmx \in \calX$ and $\calU := \calU_1 \times \cdots \times \calU_n \subseteq \mathbb{R}^n$ to denote the value space for $\utilvec(\cdot)$. 
\end{definition} 
To avoid unrealistic corner cases, we assume that $\calU = \calU(\calX)$ is a compact set. 
\begin{remark}
Many papers, for example \citet{bertsimas2011price,salem2022enabling} consider the \emph{efficiency metric} as the total utility $c (\bmx) := \sum_{i=1}^n \utility_i(\bmx)$. 
Indeed, this is central to many results derived from those papers. The assumption is natural, as maximizing $c (\bmx)$ will now correspond to maximizing the \emph{social welfare}, i.e., the total utility obtained by the stakeholders. 
However, we allow a more general efficiency metric in this paper. 
The efficiency metric could potentially be very disconnected from the utility of the stakeholders. 
For the TSP example in the introduction, the utilities of each stakeholder can be the negative of the distance they travel before reaching the hub. 
In fact, $c$ could be the negative of the entire tour, which is different from the social welfare-like metric. 
\end{remark}
Next, we motivate a metric of fairness. 
More precisely, we formally define an {\em unfairness} function that measures the degree of unfairness in the utilities obtained by the stakeholders.
}

{
\begin{definition}[Unfairness Function \citep{lodi2021fairness}]
An unfairness function $\unfair:\mathbb{R}^n\to \mathbb{R}_{\ge 0}$ is a function that satisfies the following properties.
\begin{itemize}
    \item \textbf{Perfect fairness.} $\unfair(\utility_1, \ldots, \utility_n) = 0 $ if and only if $ \utility_1 = \utility_2 = \ldots = \utility_n$.
    \item \textbf{Symmetry.}  $\unfair(\utility_1, \ldots, \utility_n) = \unfair(\pi(\utility_1, \ldots, \utility_n))$ for any permutation $\pi: [n] \to [n]$, where we use $\pi(\utility_1, \ldots, \utility_n)$ as a shorthand notation for $(\utility_{\pi(1)}, \ldots, \utility_{\pi(n)})$. 
\end{itemize}
\end{definition}
We note that quantifying fairness by an unfairness function as above is not the only way to model fairness. 
The fair allocation literature talks about other concepts like envy-freeness (each stakeholder prefers their own bundle of items compared to that of anybody else's), envy-free up to one good (each stakeholder prefers their own bundle of items compared to that of anybody else's bundle after removing at most one item from the other stakeholder's bundle), proportionality (if there are $n$ stakeholders, then each of them receives at least $1/n$-th of the maximum achievable value) and so on.
These are binary in nature, from the perspective that a solution either has this property or it does not.
In contrast, a more granular quantifier of (un)fairness as defined above is more amenable to mathematical programming formulations.

\subsection{Aggregation Functions}

The unfairness function as defined is relevant when there is a single decision-making epoch, and each stakeholder $i$ obtains a utility $\utility_i$. 
However, let us say there are $2$ rounds of decision making and stakeholder $i$ obtains a utility of $\utility_i^{(1)}$ and $\utility_i^{(2)}$ in each round, respectively. 
Then, it is not clear \emph{a priori}, if they can be aggregated in a natural way, to measure unfairness. 
It is common in the literature, 
to use the {\em averaged} utility obtained over the $T$ epochs. 
However, it is not clear if that is the most reasonable model in some applications. The issue is discussed in the following example. 
}

{
\begin{example}\label{ex:whyGenAgg}
Consider the problem of locating ambulances across a city. 
Let $x_i$ be the number of minutes it takes for the nearest ambulance to reach street (stakeholder) $i$.
The utility for street $i$ is now defined as $15-x_i$, if $x_i \le 15$ and $0$ otherwise, as an ambulance that takes a long time to reach a location, is not particularly useful. 
Suppose ambulances are to be located seven days a week. 
Let us say, the time taken for the ambulance to reach street $1$ on each of the seven days be $10,10,10,10,10,10,10$, thus making its utility $5,5,5,5,5,5,5$, and an {\em average} utility of $5$. 
In contrast, the time taken for the ambulance to reach street $2$ are $1,1,1,15,15,15,15$ on the days of the week, thus making its utility $14, 14, 14, 0, 0, 0, 0$ and an {\em average} utility of $6$, thus greater than that of street $1$.
However, it is likely that one might be more inclined to be in street 1, as opposed to street 2.
Thus, aggregating utilities obtained over multiple decision epochs by averaging need not necessarily be the best technique for every fair optimization problem.
\end{example}
}

{
    This leads to the question of how can the utilities be aggregated. In the above ambulance example,  one might think that fraction of the days when the ambulance is located so that it can reach a street within $15$ minutes is a reasonable metric of aggregation. 
    Or, one might want to reduce the variability and want to include the standard deviation of the utilities obtained over all the periods.
    In response, we define an aggregation function abstractly as follows. 
}

{
\begin{definition}[Aggregation Function]\label{def:aggFun}
$\aggFun: \left ( \bigcup\limits_{T=1}^\infty \R^T \right)  \to\R$ is an aggregation function if  it is  
\begin{enumerate}
    \item {\em Symmetric}. Namely, $\aggFun(\utility_i^{(1)}, \ldots, \utility_i^{(T)}) = \aggFun(\pi(\utility_i^{(1)}, \ldots, \utility_i^{(n)}))$ for any permutation $\pi$.  
    \item {\em Extension agnostic}. Namely, the following are equal for every $k\in\mathbb{N}$:
    \begin{itemize}
        \item $\aggFun(\utility^{(1)}_i,\ldots,\utility^{(T)}_i)$,
        \item $\aggFun(\utility^{(1)}_i,\ldots,\utility^{(T)}_i, \utility^{(1)}_i,\ldots, \utility^{(T)}_i)$,
        \item $\vdots$
        \item $\aggFun(\utility^{(1)}_i,\ldots,\utility^{(T)}_i,\underbrace{\utility^{(1)}_i,\ldots,\utility^{(T)}_i}_{k\text{ times}},\utility^{(1)}_i,\ldots,\utility^{(T)}_i) $. 
    \end{itemize}
\end{enumerate}
\end{definition} } 
{
The symmetry property of the aggregation function ensures that the aggregated value of the utilities does not vary due to the order in which these utilities are obtained. 
This enables us to talk in terms of {\em number of times} a certain utility is obtained, as opposed to the sequence of utilities themselves. 
The extension agnosticity ensures consistency in the way the utilities are aggregated. 
The utilities $\utility^{(1)}, \utility^{(2)}$ in a two-period decision making problem will have the same value as the utilities $\utility^{(1)}, \utility^{(2)}, \utility^{(1)}, \utility^{(2)}$ in a four-period decision making problem. 
Consider the extreme case of single-period decision making. Repeating the same decision for $T$ periods should not provide any additional fairness, as value is obtained primarily by the flexibility to make different decisions over time. 
The extension agnosticity property is a generalization of this property even to longer durations. 
Moreover, this property is key to the central results in the paper that provide bounds on the number of decision periods $T$ required to achieve certain levels of (un)fairness.}
\begin{remark}[Assumption of symmetry]
This assumption is no longer valid if the utility function $\bmu^{(t)}(\cdot)$ changes in each epoch $t \in [T]$, see \citet{salem2022enabling}. For instance, if agents are more sensitive to the utilities received during one time period as compared to another, such as weekends against weekdays or peak hours versus normal hours, then our framework is no longer applicable. 
We emphasize that this assumption comes {\em with} the loss of some generality. 
However, we also observe that the absence of symmetry here implies the absence of symmetry in the integer program in the natural space of variables, making the problem easier to solve with off-the-shelf integer programming technology.
\end{remark}

Next, we show explicit examples of aggregation functions.
Each of the functions mentioned below  can be represented using linear constraints with some added integer variables.

     \paragraph{The Average Function $\aggFun(\utility^{(1)},\ldots,\utility^{(T)}):=\frac{1}{T}\sum_{t=1}^T\utility^{(t)}$} The average function can be implemented by a straightforward linear function, following its definition. Actually, barring scaling and centering, the average function is the {\em only} linear aggregation function.
     \paragraph{The Minimum Function $\aggFun(\utility^{(1)},\ldots,\utility^{(T)}):=\min_{t=1}^T\utility^{(t)}$} The minimum function can be implemented by adding  a new continuous variable $y$, a set of binary variables\footnote{For ease of notation, we omit here and in the following aggregation definitions the subscript $i$ to indicate the associated stakeholder but it is easy to see that the described process has to be repeated for each $i \in [n]$.} $b_1,\ldots, b_T$ and a set of constraints $\utility^{(t)} -Mb_t \le y \le \utility^{(t)}$ for $t=1,\ldots, T$ and $\sum_t b_t= T - 1$, where $M$ is a sufficiently large positive constant. Then, variable $y$ takes the value of the aggregation function.
     \paragraph{The Maximum Function $\aggFun(\utility^{(1)},\ldots,\utility^{(T)}):=\max_{t=1}^T\utility^{(t)}$} The maximum function can be implemented by adding  a new continuous variable $y$, a set of binary variables $b_1,\ldots, b_T$ and a set of constraints $\utility^{(t)} +Mb_t \ge y \ge \utility^{(t)}$ for $t=1,\ldots, T$ and $\sum_t b_t= T - 1$, where $M$ is a sufficiently large positive constant. Then, variable $y$ takes the value of the aggregation function now.
     \paragraph{The $\rho$-th Percentile Function} We define the aggregation function as follows:
    \begin{equation*}
        \aggFun(\utility^{(1)},\ldots,\utility^{(T)}) \quad=\quad \begin{cases}
        w_{\ceil{\rho T}} & \text{ if ${\rho}{T}$ is not an integer}\\
        \frac{w_{\rho T}+w_{\rho T+1}}{2} & \text{otherwise}
        \end{cases}
    \end{equation*}
    such that $w_1,\ldots,w_T$ are a sorted permutation of $\utility^{(1)},\ldots,\utility^{(T)}$.
    Again, this function can be implemented by adding a set of new variables and constraints. Let $\utility^{(t)}$ for $t=1,\ldots,T$ be the given set of input numbers. Moreover,
    $y, w_t, b_{ts}$ are new variables added for each $t,s=1,\ldots, T$. Then, we add the following constraints: (i) $w_{t} \le w_{t+1} $, (ii) $w_t \le u^{(s)} + M(1-b_{ts})$, (iii) $w_t \ge u^{(s)} - M(1-b_{ts})$,  
    (iv) $\sum_{t=1}^T b_{ts} = 1$,  
    (v) $\sum_{t=1}^T b_{ts} = 1$, and
    (vi) $b_{ts} \in \{0,1\}$
    for $t,s=1,\dots,T$. 
    Further, if we have $\rho T$ as an integer, we add a constraint $        2v = w_{\rho T} + w_{\rho T+1} $,
    and otherwise we add $y = w_{\ceil{\rho T}}$. 
    Then, variable $y$ takes the value of the aggregation function.
     \paragraph{The Threshold Exceedance Function}Given a threshold $h \in \mathbb{R}$, we define the aggregation function as 
    \begin{equation*}
        \aggFun(\utility^{(1)},\ldots,\utility^{(T)}) \quad=\quad \sum_{t=1}^T \bm{1}\left( \utility^{(t)} \ge h \right),
    \end{equation*}
    where $1(\cdot)$ is the indicator function taking value $1$, if its argument is true and $0$ otherwise.
    This can be modeled by writing the constraints $\utility^{(t)} \le h + Mb_{t}$ and $\utility^{(t)} \ge h - \left( 1-b_{t} \right)M$, and $y = \sum_t \frac{b_{t}}{T}$. Then, $y$ takes the value of the required aggregation function. 
     \paragraph{The Mean absolute deviation Function}{We define this aggregation function as $\aggFun(\utility^{(1)},\ldots,\utility^{(T)}):= \frac{1}{T} \sum_{t=1}^T \left \vert \utility^{(t)} - \bar{\utility} \right \vert$, where $\bar{\utility}:=\frac{1}{T}\sum_{t=1}^T\utility^{(t)} $.}
    Clearly, this is a measure of variability of the utility obtained over $T$ time periods. To model this, we first compute $\bar{\utility}$ by adding a variable $\utility$ and a constraint $\utility= \sum_{t=1}^T \utility^{(t)}$. Then, we add binary variables $b_t$ and continuous variables $v_t$ for $t=1,\ldots,T$. Moreover, for each $t=1,\ldots, T$, we add the four constraints: (i) $v_t \ge  \utility^{(t)} -  \utility$, (ii) $v_t \ge -(\utility^{(t)} -   \utility)$, (iii) $v_t \le   \utility^{(t)} -   \utility + Mb_t$ and (iv) $v_t \le   \utility -   \utility^{(t)} + M(1-b_t) $. These four constraints ensure that $v_t = |  \utility^{(t)} -   \utility|$. Now, saying $y = \frac{1}{T}\sum_{t=1}^T v_t$ ensures that variable $y$ takes the value of the aggregation function.

Next, we show that certain operations among aggregation functions give an aggregation function. 
\begin{lemma}\label{thm:aggFunLinClosed}
Let $\aggFun_1, \aggFun_2, \ldots, \aggFun_k$ be $k$ aggregation functions. Then, $f(\aggFun_1, \ldots, \aggFun_k)$ is also an aggregation function for any of the following choices of $f$: 
\begin{enumerate}
    \item $f(\aggFun_1, \ldots, \aggFun_k) = \sum_{i=1}^k \alpha_i \aggFun_i$ where $\alpha_i$ are fixed real numbers. 
    \item $f(\aggFun_1, \ldots, \aggFun_k) = \max \{\aggFun_i:  1 \le i \le k\}$,
    \item $f(\aggFun_1, \ldots, \aggFun_k) = \min \{\aggFun_i:  1 \le i \le k\}$.
\end{enumerate}
Further, if $\aggFun_1, \ldots, \aggFun_k$  are mixed-integer linear representable, then so is $f(\aggFun_1, \ldots, \aggFun_k)$. 
\end{lemma}
\begin{remark}
The aggregation functions explicitly mentioned earlier, together with \cref{thm:aggFunLinClosed} imply that a variety of commonly used functions like the median, the inter-quartile range, the range etc. are all valid aggregation functions, and moreover they are all mixed-integer linear representable aggregation functions. Moreover, a suitable linear combination of, for example, the median and the range or the mean and the range could also be used as an aggregation measure. 
The motivation behind the aggregation function is that it helps in combining the utility obtained by any single stakeholder over $T$ time periods. 
If the stakeholder wants to maximize their utility in any period, the linear aggregation function (i.e., the sample mean) can be used. 
However, even if the stakeholder is interested in maximizing the average utility, while not suffering too much inconsistency, the aggregation function could be a linear combination of the sample mean and the range. 
Further, service level is generally measured in terms of percentiles, to estimate the fraction of days when high quality service is delivered or so. 
The ability to express percentiles as a mixed-integer linearly representable aggregation function directly aids us in modeling this.
\end{remark}
}
Given the above definitions, we define the \emph{fairness over time problem} as follows:
\begin{subequations}
\begin{align}
    \min_{\bm{x^{(1)}},\ldots,\bm{x^{(T)}}} ~ & ~ \unfair\left( \bmy_1,\ldots,\bmy_n \right) \label{eq:FOTA:obj}\\
    \textup{s.t.} ~&~ \bmx^{(t)} \in \calX ~~ \textup{or} ~~ (\basea) & \forall ~t \in [T] \label{eq:FOTA:feas} \\
    &\bmy_i = \aggFun \left( \bmu_i(\bmx^{(1)}), \ldots, \bmu_i(\bmx^{(T)}) \right) & \forall~i \in [n] \label{eq:FOTA:aggr}
\end{align} \label{eq:FOTA}
\end{subequations}
In problem \cref{eq:FOTA}, 
\begin{itemize}
    \item Constraints \cref{eq:FOTA:feas} ensure that the decisions made in each time period are feasible and ``sufficiently" efficient. Here, for any $\alpha \in (0,1)$, $\basea := \{\bmx \in \calX ~|~ c(\bmx) \geq \alpha \cdot \text{opt}\}$ denotes the set of $\alpha$-efficient feasible decisions with $\text{opt}$ being the optimal value of problem~\eqref{eq:dm-opt}. 
    \item Constraints \cref{eq:FOTA:aggr} aggregate for each stakeholder the (dis)utilities obtained in each time period based on the corresponding decisions.
    \item Finally, the aggregated (dis)utility vector $\bmy = (\bmy_1, \ldots, \bmy_n)$ for the $n$ stakeholders is the argument of the unfairness function \cref{eq:FOTA:obj} and the overall unfairness is minimized.
\end{itemize}
{
\citet[Example 7]{lodi2021fairness} formally demonstrates the benefit of repeated decision-making in enhancing fairness. 
This has further been noted in \citet{salem2022enabling}. 
Next, we provide estimates on what is the maximum benefit that can be obtained, and how many decision epochs are required to achieve that maximum benefit.
}

\section{The probability-equivalent reformulation} \label{sec:PER}
From a computational standpoint, the formulation in \cref{eq:FOTA} presents two fundamental issues. 
First, the feasible set and the objective function of \cref{eq:FOTA} are symmetric. 
The symmetry in the feasible set follows from the observation that any permutation of $\bmx^{(1)},\dots,\bmx^{(T)}$ is also feasible.
This is because \cref{eq:FOTA:feas} considers one $x^{(t)}$ at a time and the aggregation function used in \cref{eq:FOTA:aggr} is symmetric. 
The symmetry in the objective function follows from the fact that unfairness functions are symmetric (permutation invariant). 
Thus, all permutations of any sequence of decisions $\bmx^{(1)},\dots,\bmx^{(T)}$ are feasible and give the same objective value. 
It is well known that out-of-the-box branch-and-cut solvers for mixed-integer programs are particularly slow when handling problems with a high degree of symmetry and multiplicity of optimal solutions. 
Second, the formulation already decides the time horizon $T$ over which fairness has to be optimized. 
It is possible that for different but acceptable values of $T$, significantly less unfairness is possible. 
Conversely, a lower bound on the least possible unfairness over all values of $T$ might remain unknown. 
In light of these shortcomings, we provide what we refer to as the \emph{probability-equivalent formulation}.
\subsection{Probability-equivalent formulation}

In this section, we first define the terms and expressions required to present the probability-equivalent formulation, and then we present the formulation.

 For the rest of the paper, we use $\mathscr{D}$ to denote the set of all probability distributions with \emph{finite} support and \emph{rational} probability values. 
 We use $\bar{\mathscr{D}}$ to denote the set of all probability distributions with \emph{finite} support. 
 For any probability distribution $p_X \in  \mathscr{D} ~ \textup{or} ~ \bar{\mathscr{D}}$, we define 
$\mathscr{D}_{\supp(p_X)} ~ \textup{or} ~\bar{\mathscr{D}}_{\supp(p_X)} := \left\{ p \in 
\mathscr{D} ~ \textup{or} ~ \bar{\mathscr{D}} ~|~ \supp(p) = \supp(p_X)  \right\}
$, respectively. In other words, it is the set of probability distributions with the same support as $p_X$.

\begin{definition}[Equivalent Discrete Distribution (EDD) and Equivalent Value Support (EVS)]
    Let $\utility_i^{(1)}, \ldots, \utility_i^{(T)}$ be a finite sequence of real numbers. 
    We define the \emph{Equivalent Value Support (EVS)} $\{v_i^1,\dots, v_i^k\}$  as the set of unique numbers in the given sequence.
 Further, if $v_i^j$ is repeated $\ell_i^j$ times in the sequence for $j=1,\ldots,k$, then, we define the associated \emph{Equivalent Discrete Distribution (EDD)} as the probability distribution with probabilities, $p(v^j_i) = \frac{\ell_i^j}{T}$ for $j=1,\ldots, k$ and $p(x) = 0$ if $x \not\in \{v_i^1,\dots,v_i^k\}$.  
\end{definition}

\begin{definition}[Equivalent Enumeration Sequence (EES)]
    Let $p(x)$ be a probability mass function with finite support with rational probabilities with denominator at most $T$. i.e., for any $x$, $p(x)$ is of the form $\frac{\ell}{T}$ where $\ell$ is an integer. 
    Let the finite support be given by $\{v_1,v_2,\dots,v_k\}$ with $v_1 < v_2 < \dots<v_k$ with $p(v_i) = \frac{\ell_i}{T}$
    Then, the associated \emph{equivalent enumeration sequence (EES)} is $(w_1, w_2, \dots,w_T) := \left( \underbrace{v_1,\dots,v_1}_{\text{repeated $\ell_1$ times}}, \ldots, \underbrace{v_k,\dots,v_k}_{\text{repeated $\ell_k$ times}}  \right)$. 
\end{definition}

\begin{lemma}\label{thm:symmBreak}
    Let $\utility_i^{(1)}, \ldots, \utility_i^{(T)}$ be a finite sequence of real numbers and let $p_X$ be the EDD associated with this sequence. 
    Then the EDD associated with any permutation of $\utility_i^{(1)}, \ldots, \utility_i^{(T)}$ is also $p_X$.
\end{lemma}

{

\begin{definition}\label{defn:dist-aggFun}
\textbf{Distributional Aggregation Function.} 
We say that $\tilde{\aggFun}: \mathscr{D} \to \mathbb{R}$ is a distributional aggregation function corresponding to the aggregation function $\aggFun$ if for any finite sequence $\utility_i^{(1)}, \ldots, \utility_i^{(T)}$ and its corresponding EDD $p_X$, we have $\tilde{\aggFun}(p_X) = \aggFun\left( \utility_i^{(1)}, \ldots, \utility_i^{(T)} \right)$. 
\end{definition}

The above \cref{defn:dist-aggFun} does not say much about the existence of distributional aggregation functions. 
However, we now prove that {\em every} aggregation function $\aggFun$ has a corresponding distributional aggregation function $\tilde\aggFun$ such that $\aggFun$ evaluated at a sequence equals the value of $\tilde\aggFun$ at the EDD corresponding to the sequence. 
Conversely, the value of $\tilde\aggFun$ evaluated at any discrete distribution matches the value of $\aggFun$ evaluated at the corresponding EES. We formalize this in \cref{thm:colGenEquiv} below.
}

\begin{theorem} \label{thm:colGenEquiv}
    Let $\aggFun$ be an aggregation function.
 Then, there exists $\tilde{\aggFun}:{\mathscr{D}}\to \mathbb{R}$ with the following properties:
    \begin{enumerate}
        \item Let $p_X(x)$ be the EDD associated with the sequence $\utility_i^{(1)}, \ldots, \utility_i^{(T)}$. Then, $\aggFun\left( \utility_i^{(1)}, \ldots, \utility_i^{(T)} \right) = \tilde \aggFun (p_X)$.
        \item Let $(w_1,\ldots, w_T) $ be the EES associated with some finite probability distribution $p_X$. Then, $\aggFun(w_1,\dots,w_T) = \tilde\aggFun (p_X)$.
    \end{enumerate}
\end{theorem}
The proof of \cref{thm:colGenEquiv} is presented in \ref{app:colGenEquiv}.

Given the existence of an appropriate $\tilde\aggFun$ corresponding to any aggregation function $\aggFun$ as shown in \cref{thm:colGenEquiv}, we write the following program, which we refer to as the \emph{probabilistic equivalent (PE) formulation}.

In the formulation below, we assume that $\basea$ could be a set with potentially infinitely many elements. 
However, in any single time, we only consider a finite subset of $\basea$. 
We assume we will be able to perform {\em column generation} and {\em row generation} as needed for solving the problem. 
We present the formulation below assuming that we are considering a finite subset $\{\bmx^1,\ldots,\bmx^k\} \subseteq \basea $.

{
We use the notation $\utility^j_i$ to denote the utility obtained by stakeholder $i$ due to the decision $\bmx^j$. 
In this formulation, we want to assign a probability distribution over the set $\{\bmx^1,\ldots,\bmx^k\} $.
This translates into a probability distribution over the set of utilities $\{\utility_i^1,\ldots,\utility_i^k\}$ for stakeholder $i$. 
We denote this probability distribution as $p_X^i$ for $i=1,\ldots, n$.
Indeed, if the original distribution is on $\bmx^j$'s, then we need for any two stakeholders $i$ and $i'$, $p_X^{i}(\utility_{i}^j)=p_X^{i'}(\utility_{i'}^j)$.
}

\begin{subequations}
\begin{align}
    \min_{p_X} ~ & ~ \unfair\left( \bmy_1,\ldots,\bmy_n \right) &\textup{s.t.} \\
    \bmy_i \quad&=\quad \tilde\aggFun \left( p_X^i \right) & \forall~i \in [n] \label{eq:FOTAcolGen:agg} \\
     \utility_i^j \quad&=\quad \bmu_i\left( \bmx^j \right) & \forall~i \in [n];~j \in [k] \label{eq:FOTAcolGen:utility} \\
    \supp(p_X^i) \quad&\subseteq\quad \{\utility_i^{1},\dots, \utility_i^{k}\}
 \label{eq:FOTAcolGen:support} & \forall~i \in [n] \\
     p_X^i\left( \utility_i^{j} \right)  \quad&= \quad p_j\label{eq:FOTAcolGen:consis}& \forall~i \in [n];~j \in [k]\\
    p_j \quad&\geq\quad 0 & \forall~j \in [k] \label{eq:FOTAcolGen:positive} \\
     Tp_j \quad&\in\quad \Z_+\label{eq:FOTAcolGen:int}& \forall~j \in [k]
\end{align} \label{eq:FOTAcolGen}
In the formulation above $\bmx^1,\ldots,\bmx^k$ are elements of $\basea$. 
Since they are chosen dynamically in a column generation-related approach, they should be viewed as parameters in \cref{eq:FOTAcolGen}. 

Based on the selected decisions $\{\bmx^1,\ldots,\bmx^k\}$, the utilities $\utility_i^j$ for all $i \in [n], j \in [k]$ are fully determined via \cref{eq:FOTAcolGen:utility}. 

Here, we use variables $\{p_X^i\}_{i = 1}^n$ to denote the probability distributions of different stakeholders. The constraint \cref{eq:FOTAcolGen:support} imposes that each of the probability distributions $p_X^i$ has its corresponding finite support $\{\utility_i^1, \ldots, \utility_i^k\}$, without which \cref{eq:FOTAcolGen} will be a semi-infinite mathematical program, which we consider to be beyond the scope of this paper. 
     As discussed earlier, to reflect the fact that the probability distribution is originally on $\{\bmx^1,\ldots,\bmx^k\}$, which reflects as probability distributions on their utilities for each stakeholder, we write \cref{eq:FOTAcolGen:consis}. This enforces the previously-stated condition, that for stakeholders $i,i'$, $p_X^{i}(\utility_{i}^j)=p_X^{i'}(\utility_{i'}^j)$.

{\cref{eq:FOTAcolGen:positive} ensure that the probabilities take a non-negative value and \cref{eq:FOTAcolGen:int} ensure that the probability values all have a denominator that divides $T$. This will be used in \cref{thm:colGenequiv} to show that \cref{eq:FOTAcolGen} is indeed a reformulation of \cref{eq:FOTA}.}

Finally, aggregated utilities $\{\bmy_i\}_{i = 1}^n$ for each stakeholder are computed by constraint \cref{eq:FOTAcolGen:agg}. 

\end{subequations}

The \emph{probabilistic-equivalent} formulation in \cref{eq:FOTAcolGen} clearly overcomes the problem of symmetry. 
On closer inspection, the problem is handled by the way the aggregation function $\aggFun$ and its counterpart in \cref{eq:FOTAcolGen}, namely $\tilde\aggFun$ behave. 
As shown in \cref{thm:symmBreak}, any permutation of a sequence of solutions maps to the same probability distribution as the original sequence. 
Since $\tilde\aggFun$ acts on the probability distribution, all symmetry in the original formulation of \cref{eq:FOTA} is absent in the new formulation. 
Thus, the reformulation, in certain cases can be easier to solve than the original problem. 
We are yet to address the sensitivity of the optimal solution to choices of $T$, which we postpone to \cref{sec:probRelax}.

We would like to point out that some other fields have also applied the idea of using probabilities to express the convex combination of different decisions (\emph{strategy mixture}). In the field of stochastic optimization, \citet{delage2019dice} provides the circumstances under which random decisions (i.e., decisions selected via a random distribution) strictly outperform optimal pure (single) decisions. Later, \citet{delage2022value} gives results on the value of randomization in mixed-integer distributional robust problems. In online learning problems, \citet{castiglioni2022unifying} takes advantage of the strategy mixtures by introducing the set of probability measures on the original decision set. 

If we are able to handle the probability constraints in \cref{eq:FOTAcolGen:support,eq:FOTAcolGen:consis}, due to the lack of symmetry, it seems like \cref{eq:FOTAcolGen} might be computationally easier to solve than \cref{eq:FOTA}.
However, in our context, it is unclear, if solving \cref{eq:FOTAcolGen} provides any insights into solving \cref{eq:FOTA}. 
The following theorem provides those insights.

\begin{theorem}\label{thm:colGenequiv}
    The mathematical programs in \cref{eq:FOTA} and \cref{eq:FOTAcolGen} are equivalent in the following sense. 
    \begin{enumerate}
        \item They have the same optimal objective value. 
        \item Given a feasible solution to \cref{eq:FOTA}, a feasible solution to \cref{eq:FOTAcolGen} with the same objective value can be constructed in time bounded by a polynomial in the size of the problem data. Conversely, given a feasible solution to \cref{eq:FOTAcolGen}, a feasible solution to \cref{eq:FOTA} with the same objective value can be constructed in time bounded by a polynomial in the size of the problem data. 
    \end{enumerate}
\end{theorem}
The proof of \cref{thm:colGenequiv} is given in \ref{app:colGenequiv}.

\subsection{Probability-equivalent formulations for selected aggregation functions}\label{sec:probRelax}
In ths section, we show the equivalent probabilistic formulations of an example set of aggregation functions. 
However, we anticipate the following important observation.
\begin{remark}\label{rem:aggFunWeird}
    If an aggregation function $\aggFun$ is mixed-integer linear representable, then that does not necessarily mean that the corresponding $\tilde\aggFun$ is mixed-integer linear representable. We will, for example, show that the \emph{mean absolute deviation} aggregation function is mixed-integer linear representable, but the corresponding $\tilde\aggFun$ is not. 
    Conversely, the threshold exceedance aggregation function requires integer variables for representation. However, the corresponding $\tilde\aggFun$ can be represented exclusively by continuous variables and linear constraints in an extended space.
\end{remark}
\paragraph{Minimization Aggregation.}Let us consider the \emph{probability-equivalent} formulation of the problem, when minimization aggregation is used.
In this case, let the EDD associated with the utilities that are being aggregated be represented by $p$. In particular, let the distribution 
assign a probability of $p^j$ to $\utility^j$. 
Then, $\tilde \aggFun = \min \{\utility^j: p^j > 0\}$. However, the strict inequality constraint is generally not used while defining a mathematical program, due to the loss of closedness in the feasible set. So, we relax the condition as $\tilde \aggFun = \min \{\utility^j: p^j \ge \varepsilon\}$, where $\varepsilon$ is a sufficiently small positive constant. 
In fact, this is analogous to the ``big-M'' ideas used commonly in practice. 

Thus, we implement the above $\tilde \aggFun$ using mixed-integer linear constraints as follows:
\begin{subequations}
    \begin{align}
    r^j \quad&\ge\quad p^j &\forall~j \\
    p^j \quad & \ge \quad \varepsilon~r^j &\forall~j\\
    \sum_{j = 1}^k b^j \quad & = \quad \sum_{j = 1}^k r^j - 1 \\
    u^j r^j \quad&\le\quad y + Mb^j &\forall~j\\
    y \quad&\le\quad u^j r^j + M(1 - r^j) &\forall~j\\
    r^j, b^j \quad&\in\quad \{0,1\} &\forall~j
    \end{align}
    \label{eq:aggCGmin}
\end{subequations}
Here $\utility^j$, as before, denotes the utility obtained  due to some decision $\bmx^j$; $M$ is a sufficiently large positive number. 
The value of the aggregation gets assigned to $y$. In \cref{eq:aggCGmin}, the first two constaints ensure that the binary variable $r^j$ is equal to $1$ if and only if $p^j$ is positive (i.e., greater than or equal to $\varepsilon$). 
The next constraint ensures that if the support of the EDD has a  size $\ell$, then exactly $\ell-1$ of the $b^j$ take value $1$.
Now, consider the next two constraints. If $r^j=0$, these reduce to $y\ge 0$ (or $y+M \ge 0$) and $y \le M$, which are trivially true always. 
On the other hand, if $r^j = 1$, then the first of these two constraints say that $\utility^j \le y + Mb^j$ and the second one implies that $y \le u^j$.
From the fact that $y \le u^j$ for all $j$ such that $r^j=1$, we know that $y \le \min_j u^j$. On the other hand, we insist that there must be at least one $b^j = 0$, when $r^j = 1$. 
Observe that, this can only happen with the $j$ associated with smallest value of $\utility^j$. 
And, that the $j$ that minimizes $\utility^j$ will have the constraint equivalent to $y \ge \utility^j$. 
Thus, indeed, $y = \min_j \utility^j$ as needed. 

\paragraph{Threshold Exceedance Aggregation.}In this case, $\tilde \aggFun (p) =  P(X \ge h) $, where $X$ is a random variable that follows the distribution $p$. 
We implement this function as follows. 
\begin{subequations}
    \begin{align}
        y \quad&=\quad\sum_{j} b_{j}p_j
    \end{align}
    \label{eq:aggCGthre}
    where $b_{j}$ is a parameter that takes value $1$ if $\utility^j \ge h$ and $0$ otherwise. 
\end{subequations}

\paragraph{Mean-absolute-deviation Aggregation.}
Let the aggregation function of interest be the mean-absolute-deviation. 
In other words, let us say 
\begin{align*}
    \aggFun(\utility^{(1)}, \ldots, \utility^{(T)}) 
    := & ~ 
    \frac{1}{T} \sum_t \left \vert 
    \utility^{(t)} -  \sum_s \utility^{(s)}
    \right \vert
    \\
    \implies
    \tilde \aggFun(p)
    = & 
     \sum_j p^j \left \vert 
    \utility^{j} -  \sum_s p^s \utility^{s}
    \right \vert.
\end{align*}
Note that in this case, while $\aggFun$ is a linear function of $\utility^{(t)}$, $\tilde\aggFun$ is \emph{not} a linear function of the probabilities $p^j$. 
The following set of non-linear constraints define $\tilde\aggFun$. 
\begin{subequations}
    \begin{align}
    m \quad&=\quad \sum_{j} p^j \utility^{j} \\
    s^j \quad&\ge\quad \utility^{j} - m  &&\quad ~j\in[k]\\
    s^j \quad&\ge\quad m - \utility^{j}   &&\quad ~j\in[k]\\
    s^j \quad&\le\quad \utility^{j} - m + Mb^j  &&\quad ~j\in[k]\\
    s^j \quad&\le\quad m - \utility^{j} + M(1-b^j)  &&\quad ~ j\in[k]\\
    y \quad&=\quad \sum_{j} p^j s^j 
    \end{align}
    \label{eq:aggCGdevn}
    Clearly, in the above $m$ contains the \emph{mean} of the utility.
    The next two constraints ensure that $s^j \ge |u^j - m|$.
    And the next two constraints ensure that $s^j \le |u^j - m|$. 
    The last constraint aggregates the above and assigns the value of the required aggregation function to $y$.
\end{subequations}

We consider a problem in its entirety in the example below. 
\begin{example}\label{ex:complete}
For this example, we use $\aggFun$ as a linear combination of the minimization aggregation function and the average aggregation function. 
We use $\unfair (\bmy) = \max_i y_i - \min_i y_i$. 

Thus, the analog of \cref{eq:FOTA} is given as 
\begin{subequations}
    \begin{alignat}{40}
        \min \quad&\phantom{:~~}\quad f-g &&\quad\text{s.t.}\\
    w_1y^1_i + w_2y^2_i \quad&\ge \quad g  &&\quad \forall~i\in[n]\\
    w_1y^1_i + w_2y^2_i \quad&\le \quad f  &&\quad \forall~i\in[n]\\
    y^1_i \quad&\ge\quad u_i^{(t)} - Mb_{it}
    &&\quad \forall ~i\in[n];\; t\in [T]\\
    y^1_i \quad&\le\quad u_i^{(t)} 
    &&\quad \forall ~i\in[n];\; t\in [T] \\
    \sum_t b_{it} \quad&=\quad T-1
    &&\quad \forall ~i \in [n] \\ 
    y^2_i \quad&=\quad \frac{1}{T}\sum_{t}  u_i^{(t)}  && \quad \forall ~i \in [n] \label{eq:average}\\
    u_i^{(t)} \quad&=\quad \bmu_i (\bmx^{(t)})
    &&\quad \forall ~i\in[n];\; t\in [T] \\
    \bmx^{(t)} \quad&\in\quad \basea
    &&\quad \forall ~i\in[n];\; t\in [T] \\
    b_{it}\quad&\in\quad \{0,1\}
    &&\quad \forall ~i\in[n];\; t\in [T] 
    \end{alignat}
\end{subequations}
In the above formulation, the first two sets of constraints, along with the objective function define the unfairness function.
The next three sets of constraints define the {\em minimization} part of each stakeholder's aggregation function, while constraints \eqref{eq:average} take care of the {\em average} part of the aggregation function.
The remaining constraints ensure consistency of the formulation as in \cref{eq:FOTA}.

In this context, the \emph{probabilistic-equivalent formulation} is as follows.\footnote{ 
For convenience in writing the dual program corresponding to the linear relaxation of \eqref{eq:minPrimal}, we mention the dual variable corresponding to each constraint in parenthesis.}
\begin{subequations}
\begin{alignat}{30}
    \min_{y, b, p, r} \quad&\phantom{:~~}\quad f-g &&\quad\text{s.t.}\\
    w_1y^1_i + w_2y^2_i \quad&\ge \quad g &&\quad (\alpha_i) &&\quad \forall~i\in[n]\\
    w_1y^1_i + w_2y^2_i \quad&\le \quad f &&\quad (\beta_i) &&\quad \forall~i\in[n]\\
    u^j_i r^j \quad&\le\quad y_i^1 + Mb_i^j\qquad&&\quad (\gamma_{ij}) &&\quad \forall ~i\in[n];\; j\in [k]\\
    y^1_i \quad&\le\quad u_i^j r^j + M(1 - r^j) &&\quad (\delta_{ij})  &&\quad \forall ~i\in[n];\; j\in [k] \\
    y^2_i \quad&=\quad \sum_{j} p^j u_i^j && \quad (\upsilon_i) && \quad \forall ~i \in [n] \\
    \sum_{j = 1}^k b_i^j \quad&= \quad \sum_{j = 1}^k r^j - 1
    &&\quad (\zeta_{i})  && \quad \forall ~ i \in [n] \\
    r^j \quad&\ge\quad p^j&&\quad (\eta_j)  &&\quad \forall ~j\in [k] \\
    p^j \quad & \ge \quad \varepsilon~{ r^j}&&\quad (\theta_j) &&\quad \forall ~j\in [k] \\
    \sum _{j\in[k]} p^j \quad&=\quad 1 &&\quad (\nu) \\
    r^j, b_i^j \quad&\le\quad 1&&\quad (\lambda_j, \mu_{ij}) \qquad&&\quad \forall ~i\in[n];\; j\in [k]  \\
    Tp^j \quad&\in\quad \mathbb{Z}
    &&\quad  &&\quad \forall~ j\in [k] 
    \label{eq:minPrimal:TpjInZ}
    \\
    r^j, b_i^j \quad&\in\quad \mathbb{Z} 
    &&\quad  &&\quad \forall ~i\in[n];\; j\in [k] \label{eq:otherInteger}
\end{alignat}\label{eq:minPrimal}
\end{subequations}
\end{example}

\section{The Probabilistic Relaxation} \label{sec:ProbRel}

Clearly, \cref{eq:FOTAcolGen} is a mixed-integer program. Further, if $\tilde \aggFun$ can be written using linear and integer constraints, with potentially a few auxiliary variables, then \cref{eq:FOTAcolGen} is a mixed-integer linear program. The fact that enforcing a scaling of a variable $Tp_j$ to be an integer is not special, as one could always rewrite this as $q_j = Tp_j$ and enforce that $q_j \in \Z$.
We relax the integrality requirements on \cref{eq:FOTAcolGen} and we refer to this relaxation as the \emph{probabilistic relaxation} (PR). 
Solving the relaxed version of the problem, as we show in this subsection, gives a very interpretable solution to the problem in hand. 
In particular, this addresses the second motivation in the beginning of the section, on deciding if a carefully chosen value of $T$ could provide fairness guarantees. This is formalized by the two following theorems.

\begin{theorem}\label{thm:OptExistTh}
    The problem in \cref{eq:FOTA} can achieve an optimal objective value of $0$ for some value of $T$ (i.e., perfect fairness) if and only if the \emph{probabilistic relaxation} has an optimal objective value of $0$ and the probability distributions associated with the optimal solution have the corresponding probability masses as rational numbers. 
\end{theorem}

The proof of \cref{thm:OptExistTh} is proposed in \ref{app:OptExistTh}.

\begin{theorem} \label{thm:anyExistRat}
    Let some feasible solution (not necessarily optimal) to the \emph{probabilistic relaxation} have an objective value $\hat \unfair$ achieved by a \emph{rational} probability distribution. 
    In particular, let the distribution be $(p_1,\dots, p_k) = \left(\frac{a_1}{b_1}, \ldots, \frac{a_k}{b_k} \right)$  for some $k>0$, where $\gcd(a_j, b_j) = 1$.
    Then, \cref{eq:FOTA} has an optimal objective value equal to $\hat \unfair$ for some value of $T \le \lcm \left( b_1, \dots, b_k \right)$.
\end{theorem}

The proof of \cref{thm:anyExistRat} is presented in \ref{app:anyExistRat}.

\cref{thm:anyExistRat} establishes the strong relationship between the problem \cref{eq:FOTA} and its probabilistic relaxation provided that the distribution obtained in the probabilistic relaxation is rational. 
However, the formulation of the probabilistic relaxation does not ensure the rationality of the resulting distribution \emph{per se}. Therefore, in order to complete the story, we need to consider/discuss the case when the distribution is irrational, i.e., $p \in \bar{\mathscr{D}} \backslash \mathscr{D}$. To do that, we first introduce the following definition of distance between discrete probability distributions in $\bar{\mathscr{D}}$.
\begin{definition}\label{def:dist}
Given $p_1,p_2 \in \bar{\mathscr{D}}$, we define 
\begin{align*}
    \dist(p_1, p_2) := \begin{cases}
    \max\{|p_1(u)-p_2(u)|:u\in\supp(p_1)\} &\text{if }\supp(p_1) = \supp(p_2) \\
    \infty & \text{if }\supp(p_1) \neq \supp(p_2).
    \end{cases}
\end{align*}
Moreover, given $\bar p\in\bar{\mathscr{D}}$ and $\delta \in \R_+$, we use $\balldist(\bar p, \delta) := \{ p \in \bar{\mathscr{D}}: \dist(p, \bar{p}) < \delta \} $ to denote the neighborhood (in $\bar{\mathscr{D}}$) of $\bar{p} \in\bar{\mathscr{D}}$ with radius $\delta$. 
\end{definition}

Based on the above distance metric, we will see in \cref{lemma:continuous-aggFun} that the domains of some distributional aggregation functions $\tilde{\aggFun}$ can be extended to the set of {\em all} probability distributions with discrete support.
We will also show that, such extension can be done, so that the new function with extended support is locally Lipschitz continuous. To do that, we need the following definition.

\begin{definition} 
\label{def:real-extension}
Let $\tilde\aggFun:\mathscr{D}\to\R$ be a distributional aggregation function. 
We define the {\em real extension of $\tilde \aggFun$} as the function $\hat \aggFun : \bar{\mathscr{D}}\to \R$ such that for any probability distribution $p_X\in\mathscr{D}$ with rational probabilities, $\hat{\aggFun}(p_X) = \tilde{\aggFun}(p_X)$, and for any probability distribution $p_X \in \bar{\mathscr{D}}\setminus \mathscr{D}$, 
\begin{align}
    \hat{\aggFun}(p_X) \quad& := \quad {\liminf_{\substack{p \to p_X \\ p\in\mathscr{D}}}} \tilde{\aggFun}(p).
\end{align}
We use $\liminf$ instead of $\lim$ so that  $\hat\aggFun$ is well-defined all over $\bar {\mathscr{D}}$, and it is guaranteed to match the value of the limit when it exists.
\end{definition}

We are now ready to prove our structural lemma that states the continuous properties of typical examples of distributional aggregation functions.

\begin{lemma} \label{lemma:continuous-aggFun} 
    Let $\tilde\aggFun$ be the distributional aggregation function corresponding to the minimum or maximum aggregation, $\rho$-th percentile aggregation or the mean absolute deviation aggregation functions. Then, for any probability distribution $p_X$ with support set $\{v_1, \ldots, v_k\}$ such that $v_1 < v_2 < \dots < v_k$, and rational probability $p_X(v_i) > 0$ for all $i \in [k]$, the corresponding real extension $\hat{\aggFun}$ is $\dist(\cdot, \cdot)$-Lipschitz continuous in a ball $\balldist(p_X, \epsilon)$, where 
    \begin{enumerate}
        \item $\epsilon =p_X(v_1)$ or $p_X(v_k)$ if the aggregation function is the minimum or maximum aggregation, respectively.
        \item $\epsilon = \frac{\delta(\rho)}{k}$ if the aggregation function is the $\rho$-th percentile aggregation, with 
        \begin{align*}
            \delta(\rho) := \min \left\{ \min_{j = 1}^{k - 1} \left \{ \left| \rho - \sum_{i = 1}^{j - 1} p_X(v_i) \right|, ~  \left| \sum_{i = 1}^{j} p_X(v_i) - \rho \right| \right\}, ~ \min_{i = 1}^k p_X(v_i) \right\}  \geq 0,
        \end{align*}
        and the equality holds (i.e., $\delta(\rho) = 0$) when $\rho \in \{\sum_{i = 1}^j p_X(v_i)\}_{j = 1}^{k - 1}$. 
        \item $\epsilon = \infty$ if the aggregation is the mean absolute deviation aggregation function.
    \end{enumerate}
\end{lemma}
The proof of \cref{lemma:continuous-aggFun} is demonstrated in \ref{app:continuous-aggFun}.

\cref{lemma:continuous-aggFun} demonstrates that some distributional aggregation functions can be extended to real spaces $\bar{\mathscr{D}}$ such that they are locally Lipschitz continuous. 
Note that most of the widely used (un)fairness functions are (globally) Lipschitz continuous over the utility vector $\bmy = (\bmy_1, \ldots, \bmy_n)$, like the gap between maximum and minimum utilities $\max_{i = 1}^n \bmy_i - \min_{i = 1}^n \bmy_i$, the quadratic (un)fairness function $\bm{y}^{\top} \bm{Q} \bm{y}$ for some $\bm{Q} \succeq \bm{0}$, and so on. Thus, if the real extension of aggregation function $\hat{\aggFun}: \bar{\mathscr{D}} \to \mathbb{R}$  is (locally) Lipschitz continous with Lipschitz constant $L_{\hat{\aggFun}}$ over $\dist(\cdot, \cdot)$ (as shown in \cref{lemma:continuous-aggFun}), and the fairness function $\phi: \mathbb{R}^n \to \mathbb{R}$ is (globally) Lipschitz continuous with Lipschitz constant $L_{\phi}$ over $\dist(\cdot, \cdot)$, their composition $\phi \circ \hat{\aggFun}: \bar{\mathscr{D}} \to \mathbb{R}$ is (locally) Lipschitz continuous with Lipschitz constant $L_{\phi} \cdot L_{\hat{\aggFun}}$ over $\dist(\cdot, \cdot)$, i.e., the product of the Lipschitz constants.

We are now ready to state and prove our main result, i.e., that one can construct a sequence of decisions for the original fairness over time problem that achieve a near-optimal objective value based on the probability vector obtained by solving the probabilistic relaxation. For such result, given in \cref{thm:anyExistLips} below, we assume that the (composite) objective function is (locally) Lipschitz continuous over $\bar{\mathscr{D}}$, which holds with typical (un)fairness and aggregation functions.

\begin{theorem} \label{thm:anyExistLips}
    Consider an unfairness function $\unfair$ and a distributional aggregation function $\tilde \aggFun$. Let the real extension of $\tilde \aggFun$ be denoted by $\hat\aggFun$.
    Suppose the composition of the unfairness function and the aggregation function, i.e., $\unfair(\hat \aggFun(p_X^1),\ldots, \hat \aggFun(p_X^n) )$ is Lipschitz continuous with Lipschitz constant $L$. 
    Then, given a solution to the probabilistic relaxation with objective value $\hat \unfair$, a solution to \cref{eq:FOTA} with objective value at most $\hat\unfair + \epsilon$ can be constructed in time bounded by a polynomial in the problem data for any given value of $T \geq \ceil{\frac{L}{\epsilon}}$. 
\end{theorem}
The proof of \cref{thm:anyExistLips} is given in \ref{app:anyExistLips}.

We end the section by reprising \cref{ex:complete} and presenting explicitly the associated probabilistic relaxation.

\begin{example}[Continuation of \cref{ex:complete}.]
One can obtain the probability relaxation of the problem considered in \cref{ex:complete} by removing constraint \cref{eq:minPrimal:TpjInZ} from \cref{eq:minPrimal}. 
While this removes some of the integrality requirements from the probability-equivalent formulation, this is {\em still} a mixed-integer linear program.
However, this does not suffer anymore from symmetry,  allows us a choice of $T$, and can be solved using branch-and-price algorithms.

To compute the column generation for the branch-and-price methods, we present the dual program corresponding to the linear relaxation of \cref{eq:minPrimal}, i.e., that obtained by removing also \eqref{eq:otherInteger}, as follows:
\begin{subequations}
\begin{alignat}{30}
 \max 
  \quad\phantom{:}\quad
    -\sum_j \lambda_j - \sum_{i,j}( \mu_{ij} &+ 
     M \delta_{ij} ) 
    +\sum_i 
    \zeta_i
    + \nu &&\quad \text{s.t.} \\
  \sum_i \beta_i \quad&\le\quad 1 &&(f) 
  \label{eq:minDual:f}
  \\
    \sum_i \alpha_i \quad&\ge\quad 1 &&(g) 
    \label{eq:minDual:g}
    \\
  w_1(\alpha_i - \beta_i) + \sum_j (\gamma_{ij} - \delta_{ij})
    \quad&\le\quad 0 &&(y^1_i) &\forall~i \in [n] 
    \label{eq:minDual:yOne}
    \\
  w_2(\alpha_i - \beta_i) + \upsilon_i 
    \quad&\le\quad 0 &&(y^2_i) &\forall~i \in [n] 
    \label{eq:minDual:yTwo}
    \\
  -\sum_i \utility_i^j \upsilon_i - \eta_j + \theta_j + \nu 
    \quad&\le\quad 0 && (p_j)
    &\forall~j \in [k]
    \label{eq:minDual:pj}
    \\
  \sum_i \utility_i^j (\delta_{ij}-\gamma_{ij}) 
   + \sum_i \zeta_i + \eta_j - \varepsilon \theta_j -\lambda_j 
    \quad&\le\quad 0 && (r_j)
    &\forall~j \in [k]
    \label{eq:minDual:rj}
    \\
  M\gamma_{ij} - \zeta_i - \mu_{ij}
  \quad&\le\quad 0 &&(b_{ij})
    &\forall~ i \in [n];~
    j \in [k] 
    \label{eq:minDual:bij}
\end{alignat}\label{eq:minDual}
\end{subequations}
Here, observe that  $\lambda_j$ and $\mu_{ij}$  appear with a negative sign in the objective and the only constraints where they appear are \cref{eq:minDual:rj} and \cref{eq:minDual:bij}, respectively. 
Similarly, $\nu$ appears with a positive sign in the objective and it appears only in \cref{eq:minDual:pj}.
Thus, these variables and constraints \cref{eq:minDual:rj,eq:minDual:bij,eq:minDual:pj}  can be removed by the following substitution:
\begin{subequations}
    \begin{align}
        \mu_{ij} \quad&=\quad 
        \max \{0, 
        M\gamma_{ij} - \zeta_i \} 
        \\
        \lambda_j \quad&=\quad \max  \left\{ 0, 
        \sum_i \utility_i^j (\delta_{ij}-\gamma_{ij}) 
   + \sum_i \zeta_i + \eta_j - \varepsilon \theta_j
        \right\}
        \\
        \nu \quad&=\quad \min_j \left\{\sum_i \utility_i^j \upsilon_i + \eta_j - \theta_j\right\}
    \end{align}\label{eq:minSubs}
\end{subequations}
More importantly, these are the {\em only} constraints that vary with respect to the choice of $k$.

Given this dual, one can solve the linear programming relaxation of \cref{eq:minPrimal} with a value $k \ll |\basea|$.
This is equivalent to using only a subset of $\basea$ as feasible decisions. 
From the fact that we are allowing only a subset of strategies, this is a restriction. 
However, given the formulation, extending the set of strategies adds more constraints to the model. 
Thus, simultaneous column and row generation is required to solve this problem. 
Solving this restriction will give a suboptimal solution to the complete \cref{eq:minPrimal} where $k = |\basea|$. 
Suboptimality in the primal corresponds to infeasibility in the dual \cref{eq:minDual}, which leads to the classical column generation scheme in which, for example, one can identify the most violated constraint in \cref{eq:minDual} after the substitutions in \cref{eq:minSubs}.
This directly leads to identifying the variable that $\bmx\in \basea$ that enters the basis.

\end{example}

Following the steps in the above example, the simultaneous row and column generation scheme can be adapted for most of the aggregation and (un)fairness functions of interest. Such adaptation is direct if those functions can be written in mixed-integer linear terms. Otherwise, an extension could be significantly more complex.

\section{Discussion}\label{sec:Disc}
In this study, we provided a general framework for modeling and solving issues of fair decision making when multiple decision rounds are either natural (repeated decision making) or at least conceivable. 
Specifically, we have expanded the formulation techniques in the literature and provided a conceptualization for aggregating the utilities acquired by stakeholders over multiple time periods. 
While we have assumed that all stakeholders use the same aggregation function $\aggFun$, i.e., they aggregate the utilities in the same way, the findings still hold if we assume that each individual stakeholder $i$ uses a different aggregation functions $\aggFun_i$. 
In addition, we provided a generalized formulation (the probability-equivalent formulation) that may be solved by employing row and column generation-based strategies. 
This is especially pertinent given that the problem in its natural form is susceptible to a high degree of symmetry. 
In addition, we demonstrated that solving a particular relaxation of the probability-equivalent reformulation is highly beneficial. First, it offers insights into the number of rounds of decision making required to attain high levels of fairness. Second, the probability-equivalent approach isolates the problem of recognising fair decisions from the challenge of optimizing over the set of efficient feasible solutions $\basea$. 
For instance, if one is interested in traveling salesman tours that are fair over time, according to some measures of utility and fairness, one can generate tours as the \emph{pricing problem} associated with column generation, while retaining only the fair decision-making component in the master problem. This is in contrast to the natural approach that requires to concurrently identify TSP tours and optimizing fairness (while dealing with symmetry).

\paragraph{Limitations and future work}The lack of Pareto optimality is among the shortcomings of this approach. 
A solution $\bmx\in\basea$ is said to be Pareto-dominated, if there exists $\bmx'\in\basea$ such that $\utilvec_i(\bmx') \ge \utilvec_i(\bmx)$ for each stakeholders $i$ with strict inequality holding for at least one stakeholder. 
A solution that is Pareto-optimal is one that is not Pareto-dominated. 
Pareto-dominated solutions are not recommended since it is possible to increase the utility of certain stakeholders without harming others. 
This is a crucial requirement in fair allocation issues because, as discussed in the introduction, all stakeholders may be assigned nothing, providing everyone a utility of $0$ and making the solution fully fair. 
While it is straightforward to avoid such extremes with our framework, since we can simply impose an efficiency threshold $\alpha$ for a solution to be feasible, i.e., to belong to $\basea$, our approach does not prohibit explicitly the selection of Pareto-dominated solutions. One potential way of dealing with Pareto dominance is to redefine the set $\basea$ by removing Pareto-dominated elements, which, however, does not look trivial from an operational standpoint, so  
an abstract way to directly address this issue is definitely an interesting direction for novel research.

This paper's findings may be expanded upon in several ways. One can have an {\em online} version of this decision-making problem, where utility metrics or feasible sets can (slowly) change over time. 
Such questions are being already considered for a {\em fair resource allocation} problem \citep{bampis2018fair}. 
In the context of general decision-making problems, an online version might eliminate or at least mitigate the need of aggregation techniques, thus leading to a potential different interpretation of the reformulation schemes proposed in this paper. More interestingly, utility metrics or feasible sets at an upcoming time-period could even depend upon the decisions made in the current time period. Identifying fair decisions in these contexts could be interesting.

Another direction is to consider specific families of the set $\basea$.
For example, $\basea$ could be the set of ($\alpha$-efficient) matchings or TSP tours in a graph. It would be useful to understand what additional results hold as a consequence of our special knowledge of $\basea$.

\section*{Acknowledgements}
This work has been partially conducted when the authors were part of the Canada Excellence Research Chair (CERC) in ``Data Science for Real-time Decision-Making" at Polytechnique Montr\'eal. We warmly thank the CERC support as well as the positive influence of all CERC members. GW was supported by the Singapore MOE under AcRF Tier-1 grant 22-5539-A0001. 

\nocite{alkema2020stripTSP}
\bibliographystyle{plainnat}
\bibliography{reference.bib}

\newpage
\appendix
\noindent {\LARGE \bf Appendix} 

\section{Proofs in \cref{sec:PER}} 

\subsection{Proof of \cref{thm:colGenEquiv}} \label{app:colGenEquiv}

\begin{proof}
{ 
To show that we can construct $\aggFun$ such that for any given $\utility_i^{(1)}, \ldots, \utility_i^{(T)} $ and its corresponding EDD, $p_X$, $\aggFun(\utility_i^{(1)}, \ldots, \utility_i^{(T)} ) = \tilde\aggFun(p_X)$, all that we have to show is that different sequences that have the same EDD should also evaluate to the same value of $\aggFun$. 
Whatever this value is, that can be defined as the value of $\tilde \aggFun(p_X)$, providing the construction.
In other words, it is sufficient to prove that if $\utility_i^{(1)}, \ldots, \utility_i^{(T_i)} $ and $\utility_{i'}^{(1)}, \ldots, \utility_{i'}^{(T_{i'})} $ are two sequences with the same EDD, then $\aggFun(\utility_i^{(1)}, \ldots, \utility_i^{(T_i)} ) = \aggFun(\utility_{i'}^{(1)}, \ldots, \utility_{i'}^{(T_{i'})} )$.
}

{
Now, for these two finite sequences (not necessarily of the same length), if the corresponding EDDs are the same, then the sequences must have the same EVS, as the EVS is nothing but the support of the EDD.
Let the EVS be $\{v^1,\ldots,v^k\}$. 
Let $v^j$ be repeated $\ell^j_i$ times in the first sequence and $\ell^j_{i'}$ times in the second sequence. 
However, since the EDD is the same for both the sequences, we necessarily require that 
$p_X(v^j) = \frac{\ell^j_i}{T_i} = \frac{\ell^j_{i'}}{T_{i'}}$ for all $j=1,\ldots,k$. 
In other words, $\ell^j_iT_{i'} = \ell^j_{i'}T_i$.
}

{
Now, consider the sequence $S_i$, where $\utility_i^{(1)}, \ldots, \utility_i^{(T_i)}$ is repeated $T_j$ times and the sequence $S_{i'}$ where $\utility_{i'}^{(1)}, \ldots, \utility_{i'}^{(T_{i'})}$ is repeated $T_i$ times. 
Thus, $S_i$ has each $v^j$ repeated $\ell^j_i T_{i'}$ times, while $S_{i'}$  has each $v^j$ repeated $\ell^j_{i'} T_i$ times for $j=1,\ldots,k$. 
However, we have already proved that $\ell^j_iT_{i'} = \ell^j_{i'}T_i$.
Thus, $S_i$ and $S_{i'}$ are only permutations of each other, and by the symmetry property of aggregation functions, they evaluate to the same value, i.e., $\aggFun(S_i) = \aggFun(S_{i'})$.
}

{
Finally, due to the extension agnosticity property of aggregation functions,
$\aggFun(\utility_i^{(1)}, \ldots, \utility_i^{(T_i)}) = \aggFun(S_i)$ and 
$\aggFun(\utility_{i'}^{(1)}, \ldots, \utility_{i'}^{(T_{i'})}) = \aggFun(S_{i'})$.
However, $\aggFun(S_i)=\aggFun(S_{i'})$, indicating that 
$\aggFun(\utility_i^{(1)}, \ldots, \utility_i^{(T_i)} ) = \aggFun(\utility_{i'}^{(1)}, \ldots, \utility_{i'}^{(T_{i'})} )$.
}

{
Property 2 holds from the fact that if $(w_1,\ldots,w_T)$ is the EES associated with $p_X$, then the EDD associated with $(w_1,\ldots,w_T)$ is indeed $p_X$, and then applying the first property. 
}
\end{proof}

\subsection{Proof of \cref{thm:colGenequiv}} \label{app:colGenequiv}

\begin{proof} 
First, we prove part 2, which then implies part 1. Let $\{\bmx^{(t)}\}_{t = 1}^T$ be a feasible sequence of  decisions for \cref{eq:FOTA} giving any objective function value $\phi$.

Let us say $(\bmx_1,\ldots,\bmx_\kappa)$ be the {\em unique} decisions, where each $\bmx_j\in\basea$ is repeated $q_j$ times.
Then, $p^i_X$ in \cref{eq:FOTAcolGen} can be chosen as the 
EDD corresponding to the sequence $\utility_i(\bmx^{(1)}), \ldots, \utility_i(\bmx^{(T)})$.
Clearly, $p_X^i$ has finite support as required.
Moreover, since the number of times $\utility_i(\bmx_j)$ is repeated is $q_j$ times, which is independent of $i$, $p^i_X(\utility_i(\bmx_j)) = \frac{q_j}{T} = p_j$ is invariant for all values of $i$, thus satisfying \cref{eq:FOTAcolGen:consis}.
Moreover, since each $p_j$ is defined as $\frac{q_j}{T}$, $Tp_j = q_j \in \Z$, implying the feasibility for \cref{eq:FOTAcolGen}.

We have
\begin{align*}
    \bmy_i^{\cref{eq:FOTA}} = & ~ \aggFun(\utility_i(\bmx^{(1)}), \ldots, \utility_i(\bmx^{(T)})) \\
    = & ~ \aggFun(\underbrace{\utility_i(\bmx_1), \ldots, \utility_i(\bmx_1)}_{q_1 ~ \textup{times}}, \ldots, \underbrace{\utility_i(\bmx_{\kappa}), \ldots, \utility_i(\bmx_{\kappa})}_{q_{\kappa} ~ \textup{times}}) \\
    = & ~ \tilde{\aggFun}(p_X^i) = \bmy_i^{\cref{eq:FOTAcolGen}}.
\end{align*}
Here, the first equality holds due to the symmetric property of aggregation functions and the second equality holds due to \cref{thm:colGenEquiv}. Since the values of $\bmy_i$ are the same for both formulations, so will the objective that only depends on them. 

Conversely, given a feasible solution to \cref{eq:FOTAcolGen}, observe that $Tp_j$ will be an integer, which we call $q_j$. 
Now, construct $\bmx^{(1)},\ldots,\bmx^{(T)}$ by repeating each $\bmx_j$, $q_j$ times. As before, by \cref{thm:colGenEquiv}, we obtain the same values of $\bmy_i$'s and hence the same objective function value. 

Part 1 follows as a corollary to part 2.
\end{proof}

\section{Proofs in \cref{sec:ProbRel}} \label{app:ProbRel}

\subsection{Proof of \cref{thm:OptExistTh}} \label{app:OptExistTh}
\begin{proof}
 If the \emph{probabilistic relaxation} has an optimal objective value of $0$ and if the probability distributions associated with the optimal solution have the corresponding probability masses as rational numbers, one can construct a feasible solution of \cref{eq:FOTA} with $T = \lcm(b_1, \ldots, b_M)$ time slots of objective value $0$ using \cref{thm:colGenequiv}.  
 
 On the other hand, assume that perfect fairness can be achieved within a finite time period $T$. Let $q_1, \ldots, q_k$ be the number of times that the candidate solution $\bmx_i$ is chosen during $T$ time periods. Based on the property of distributional aggregation functions, we have $0 = \aggFun(\utility_i^{(1)}, \ldots, \utility_i^{(T)}) = \tilde{\aggFun}(p_X)$, where 
    \begin{align}
        p_X = \left( \frac{q_1}{T}, \ldots, \frac{q_k}{T} \right) \in \mathbb{Q}^k
    \end{align}
is the corresponding EDD. \qedhere
\end{proof}

\subsection{Proof of \cref{thm:anyExistRat}} \label{app:anyExistRat}

\begin{proof}
 Define $\hat{T} := \lcm(b_1, \ldots, b_k)$. Represent the distribution $\hat{\calD}$ with its probability vector over the support set 
\begin{align}
    \hat{p} = \left( \frac{\frac{\hat{T}a_1}{b_1}}{\hat{T}}, \ldots, \frac{\frac{\hat{T}a_k}{b_k}}{\hat{T}} \right), 
\end{align}
where $\frac{\hat{T}a_1}{b_1}, \ldots, \frac{\hat{T}a_k}{b_k}$ are non-negative integers with $\frac{\hat{T}a_1}{b_1} + \cdots + \frac{\hat{T}a_k}{b_k} = \hat{T}$. Now, consider the sequence of utilities 
\begin{align}
\label{eq:sequence}
    \bigg\{ \underbrace{\utility_i^{(1)}, \ldots, \utility_i^{(1)}}_{\text{$\hat{T}a_1/b_1$ times}}, \underbrace{\utility_i^{(2)}, \ldots, \utility_i^{(2)}}_{\text{$\hat{T}a_2/b_2$ times}}, \ldots, \underbrace{\utility_i^{(k)}, \ldots, \utility_i^{(k)}}_{\text{$\hat{T}a_k/b_k$ times}} \bigg\}. 
\end{align}
Due to the property of EES and distributional aggregation functions,  sequence \eqref{eq:sequence} has the EDD $(p_1, \ldots, p_k)$ that achieves the same objective value $\hat{\phi}$. 
\end{proof}

\subsection{Proof of \cref{lemma:continuous-aggFun}} \label{app:continuous-aggFun}

\begin{proof}
We prove the lemma for each of the considered aggregation functions below.

     \paragraph{Minimum and Maximum aggregation functions} First, let us consider the minimum aggregation function. We have $\tilde{\aggFun} (p_X) = \min_{i = 1}^k v_i = v_1$ with probability $p_X(v_1) > 0$. Now consider any \emph{rational} probability vector $\tilde{p}_X \in \mathscr{D} \cap \balldist(p_X, p_X(v_1))$, we have value $v_1$ in the support of $\tilde{p}_X$ by $\tilde{p}_X(v_1) > p_X(v_1) - p_X(v_1) = 0$, and thus
        \begin{align*}
            & v_1 = \tilde{\aggFun} (p_X) = \tilde{\aggFun} (\tilde{p}_X), & \forall ~ \tilde{p}_X \in \mathscr{D} \cap \balldist(p_X, p_X(v_1)).  
        \end{align*}
    Since $\tilde{\aggFun}$ is constant in $\mathscr{D} \cap \balldist(p_X, p_X(v_1))$, the real extension $\hat{\aggFun}$ satisfies 
    \begin{align*}
        & \hat{\aggFun} (\hat{p}_X) := \liminf_{\substack{p \to \hat{p}_X \\ p\in\mathscr{D} \cap \balldist(p_X, p_X(v_1)) }} \tilde{\aggFun} (p) = \lim_{\substack{p \to \hat{p}_X \\ p\in\mathscr{D} \cap \balldist(p_X, p_X(v_1)) }} \tilde{\aggFun} (p) = v_1, & \forall ~ \hat{p}_X \in \balldist(p_X, p_X(v_1)). 
    \end{align*}
    Therefore, we show that the minimum aggregation function is locally Lipschitz continuous (constant is a special case of Lipschitz continuous) over the set $\balldist(p_X, p_X(v_1))$. 

    An analogous argument shows the property for the maximum aggregation function.
 
 \paragraph{$\rho$-th percentile aggregation function}
Consider any $\rho \in [0,1]$ and choose $j$ such that  
        \begin{align*}
          p_X(v_1) + \cdots + p_X(v_{j - 1}) \leq \rho \leq p_X(v_1) + \cdots + p_X(v_{j - 1}) + p_X(v_{j}).
        \end{align*} 
        When $\rho \neq \sum_{i = 1}^j p_X(v_i)$ for all $j = 1, \ldots, k - 1$ (fewer than $k$ points in total), define
        \begin{align*}
            \delta(\rho) := \min\left\{\rho - \sum_{i = 1}^{j - 1} p_X(v_i), ~ \sum_{i = 1}^j p_X(v_i) - \rho, ~ \underbrace{\min_{i = 1}^k p_X(v_i)}_{\textup{may not need  due to $\infty$ distance of dist-norm}} \right\} > 0. 
        \end{align*}
        Now, for any $\rho \in [0, 1] \backslash \left\{\sum_{i = 1}^j p_X(v_i)\right\}_{j = 1}^{k - 1}$, any \emph{rational} probability vector $ \tilde{p}_X \in \mathscr{D} \cap \balldist(p_X, \delta(\rho) / k)$, {nonzero components in $\tilde{p}_X$ satisfies $0 < \tilde{p}_X(v_i) \in [ p_X(v_i) - \delta(\rho) / k, ~ p_X(v_i) + \delta(\rho) / k] $ for all $i \in [k]$ due to the definition of $\dist(\cdot, \cdot)$. Then adding up the first $j - 1$ or $j$ probability values implies} 
        \begin{align*}
            & \tilde{p}_X (v_1) + \cdots + \tilde{p}_X (v_{j - 1}) \leq p_X(v_1) + \cdots + p_X(v_{j - 1}) + \delta(\rho) \frac{j - 1}{k} \leq \rho, \\
            & \tilde{p}_X (v_1) + \cdots + \tilde{p}_X (v_{j}) \geq  p_X(v_1) + \cdots + p_X(v_{j - 1}) + p_X(v_{j}) - \delta(\rho) \frac{j}{k} > \rho.
        \end{align*}
        Therefore, the $\rho$-th percentile distributional aggregation function satisfies 
        \begin{align*}
            & \tilde{\aggFun} (p_X) = \tilde{\aggFun} (\tilde{p}_X) = v_j, & \forall ~ \tilde{p}_X \in \mathscr{D} \cap \balldist(p_X, \delta(\rho) /k). 
        \end{align*}
        Then, the real extension $\hat{\aggFun}$ ensures:
        \begin{align*}
            & \hat{\aggFun} (\hat{p}_X) := \liminf_{\substack{p \to \hat{p}_X \\ p\in\mathscr{D} \cap \balldist(p_X, \delta(\rho) / k) }} \tilde{\aggFun} (p) = \lim_{\substack{p \to \hat{p}_X \\ p\in\mathscr{D} \cap \balldist(p_X, \delta(\rho) / k) }} \tilde{\aggFun} (p) = v_j, & \forall ~ \hat{p}_X \in  \balldist(p_X, \delta(\rho) /k),
        \end{align*}
        which is also Lipschitz continuous over the ball $\balldist(p_X, \delta(\rho) /k)$. 
 { \paragraph{The mean absolute deviation function} 
Consider two \emph{rational} probability distributions $p_X, p_X' \in \mathscr{D}$ with $\dist(p_X, p_X') < \infty$.
In other words, $p_X$ and $p_X'$ have the same support, $v_1, \ldots, v_k$.
First we show that the means of $p_X$ and $p_X'$ are not too far away.
Let $\bar{v}_X$ and $\bar{v}_X'$ be the means of $p_X $ and $p_X'$ respectively. Then,
        \begin{align*}
            |\bar{v}_X - \bar{v}_X'| = \left\vert \sum_{j = 1}^k (p_X(j) - p_X'(j)) v_j \right\vert \leq \sum_{j = 1}^k |p_X(v_j) - p_X'(v_j)| |v_j| \leq \dist(p_X, p_X') \|v\|_1.
        \end{align*}
        Hence, the difference between $p_X, p_X' \in \mathscr{D}$ satisfies 
  {
  \begin{align*}
      \left\vert \tilde{\aggFun}(p_X) - \tilde{\aggFun}(p_X') \right\vert = & ~ \left\vert \sum_{j = 1}^k p_X(v_j) |v_j - \bar{v}_X| - \sum_{j = 1}^k p_X'(v_j) |v_j- \bar{v}_X,| \right\vert \\
      = & ~ \left\vert \sum_{j = 1}^k \left( p_X(v_j) |v_j - \bar{v}_X| - p_X'(v_j) |v_j - \bar{v}_X'| \right) \right\vert \\
     \leq & ~ \sum_{j = 1}^k 
    \left\vert \left( p_X(v_j) |v_j - \bar{v}_X| - p_X'(v_j) |v_j - \bar{v}_X'| \right) \right\vert \\
    \leq & ~ \sum_{j = 1}^k \left\vert  p_X(v_j) |v_j - \bar{v}_X| - p_X'(v_j) |v_j - \bar{v}_X| \right\vert + p_X'(v_j) |\bar{v}_X - \bar{v}_X'| \\
    = & ~ \sum_{j = 1}^k \left\vert  (p_X(v_j) - p_X'(v_j))|v_j - \bar{v}_X| \right\vert + p_X'(v_j) |\bar{v}_X - \bar{v}_X'| \\
    \leq & ~ \dist(p_X, p_X') \|v - \bar{v}_X \bm{1}\|_1 + |\bar{v}_X - \bar{v}_X'|.
  \end{align*}
  Here, the final inequality holds by using Cauchy-Schwarz inequality with $\ell_1$-norm, $\ell_{\infty}$-norm for the first term summation, i.e., 
  \begin{align*}
      \sum_{j = 1}^k \left\vert  (p_X(v_j) - p_X'(v_j))|v_j - \bar{v}_X| \right\vert = & ~ \sum_{j = 1}^k \left\vert  (p_X(v_j) - p_X'(v_j))\right\vert \cdot |v_j - \bar{v}_X| \\
      \leq & ~ \max_{j = 1}^k\{|p_X(v_j) - p_X'(v_j)|\} \cdot \|v - \bar{v}_X \bm{1}\|_{1} \\
      = & ~ \dist(p_X, p_X') \cdot \|v - \bar{v}_X \bm{1}\|_1
  \end{align*}
  with $\max_{j = 1}^k\{|p_X(v_j) - p_X'(v_j)|\} = \dist(p_X, p_X')$ when $p_X, p_X'$ share the same support $\{v_1, \ldots, v_k\}$. Inserting the result $|\bar{v}_X - \bar{v}_X'| \leq \dist(p_X, p_X') \|v\|_1$ into the previous inequality implies 
  \begin{align*}
       \left\vert \tilde{\aggFun}(p_X) - \tilde{\aggFun}(p_X') \right\vert \leq & ~ \dist(p_X, p_X') \|v - \bar{v}_X \bm{1}\|_1 + |\bar{v}_X - \bar{v}_X'|\\
       \leq & ~ \dist(p_X, p_X') \left( \|v - \bar{v}_X \bm{1}\|_1 + \|v\|_1 \right) \\
            \leq & ~ \dist(p_X, p_X') \left( 2 \|v\|_1 + \max_{j = 1}^k \|v_j \bm{1}\|_1 \right)
  \end{align*}
        Now, the real extension $\hat{\aggFun}$ satisfies that: for any $\hat{p}_X \in \bar{\mathscr{D}}$, we have
        \begin{align*}
            \hat{\aggFun}(\hat{p}_X) := \liminf_{\substack{p \to \hat{p}_X \\ p\in\mathscr{D} }} \tilde{\aggFun} (p) = \lim_{\substack{p \to \hat{p}_X \\ p\in\mathscr{D} }} \tilde{\aggFun} (p), 
        \end{align*}
        where we claim the right-hand side limit exists since every convergence sequence $\{p_n\} \rightarrow \hat{p}_X$ in $\mathscr{D}$ forms a Cauchy sequence. The such claim holds since the difference between the two rational distributions in $\mathscr{D}$ is upper bounded with respect to the product between some fixed constant $\|v - \utility_i \bm{1}\|_1 + \|v\|_1$ and the distance $\dist(p_n, \hat{p}_X) \rightarrow 0$ as $n \rightarrow \infty$ for any sequence $\{p_n\}$. Therefore, similar to $\tilde{\aggFun}$, for any two points $\hat{p}_X, \hat{p}_X' \in \bar{\mathscr{D}}$, we have $|\hat{\aggFun}(\hat{p}_X) - \hat{\aggFun}(\hat{p}_X')| \leq \dist(\hat{p}_X, \hat{p}_X')(2 \|v\|_1 + \max_{j = 1}^k \|v_j \bm{1}\|_1 )$, i.e., Lipschitz continuous. 
        \qedhere
}}
\end{proof}

\subsection{Proof of \cref{thm:anyExistLips}} \label{app:anyExistLips}

\begin{proof}
    Let $\{p_X^i\}_{i = 1}^n$ be the probability distributions with respect to the objective value $\hat{\phi}$ and support set $\{v_1^i, \ldots, v_k^i\}$ for all $i \in [n]$. Notice that for any  $i \neq i' \in [n]$, we always have $p_X^i(v_j^i) = p_X^{i'}(v_j^{i'})$ as required in \cref{eq:FOTAcolGen}. To be concise, we use $p_X(v_j) = p_j$ in Problem~\cref{eq:FOTAcolGen} to denote the value of probability with respect to $j$-th component by ignoring the superscript. For any index $j \in [k]$, let $\tilde{T}_j := \lfloor T \cdot p_X(v_j) \rfloor$, and let $\tilde{T} := \sum_{i = 1}^k \tilde{T}_j \leq T$. Now, for all $j \in [k]$, define 
\begin{align}
    T_j := \left\{ 
    \begin{array}{lll}
        \tilde{T}_j + 1 & \text{if } j = 1, \ldots, T - \tilde{T} \\
        \tilde{T}_j & \text{if } j = T - \tilde{T} + 1, \ldots, k \\
    \end{array}
    \right. . 
\end{align}
Notice that the above setting ensures that $\sum_{j = 1}^k T_j = T$. Now, construct the probability distribution $p_{\mathcal{S}}^i$ using $\{T_j\}_{j = 1}^k$ proposed above such that 
\begin{align}
    & p_{\mathcal{S}}^i(v_j^i) = \frac{T_j}{T}, & \forall ~ j \in [k]
\end{align}
with the same support as $p_X^i$ for all $i \in [n]$. Then, we need to verify the distance between $p_X^i$ and $p_{\mathcal{S}}^i$ for each $i \in [n]$. Since, 
no matter $i \in [n]$, for each $j \in [k]$, 
\begin{align}
    |p_X^i(v_j^i) - p_{\mathcal{S}}^i(v_j^i)| = \left\{ 
    \begin{array}{lll}
        \left| \frac{\lfloor T \cdot p_X(v_j) \rfloor + 1 - T \cdot p_X(v_j)}{T} \right| \leq \frac{1}{T} & \text{if } j = 1, \ldots, T - \tilde{T}  \\
        \left| \frac{\lfloor T \cdot p_X(v_j) \rfloor - T \cdot p_X(v_j)}{T} \right| \leq \frac{1}{T} & \text{if } j = T - \tilde{T} + 1, \ldots, k 
    \end{array}
    \right.,
\end{align}
then we have $\dist(p_X^i, p_{\mathcal{S}}^i) = \max_{j = 1}^k |p_X^i(v_j^i) - p_{\mathcal{S}}^i(v_j^i)|  \leq \frac{1}{T}$ for all $T \in \mathbb{Z}_+$. Now, in order to satisfy $\unfair_T - \hat{\unfair} \leq \epsilon$, let $T \geq \ceil{\frac{L}{\epsilon}}$, then we have 
\begin{align*}
    \dist(p_X^i, p_{\mathcal{S}}^i) \leq \frac{1}{T} \leq \frac{\epsilon}{L}, 
\end{align*}
which implies that $\unfair_T - \hat{\unfair} \leq L \max_{i = 1}^n \{\dist(p_X^i, p_{\mathcal{S}}^i)\} \leq L \frac{1}{T} \leq \epsilon$. 
\end{proof}

\end{document}